Memorial for Dr. Erik Bollt

Scott Bollt
Mario di Bernardo
Jeremie Fish
Daniel J. Gauthier
Ying-Cheng Lai
Mason A. Porter
Jie Sun

## Introduction

Applied mathematician Erik M. Bollt, who produced a diverse body of research in applied dynamical systems and complex systems, died suddenly on 7 December 2025 while cross-country skiing. Erik was born on 24 April 1967 in Washington, DC, and grew up in Bethesda, MD, USA. He earned a B.S. degree in Applied Mathematics from the University of California at Berkeley in 1990, an M.S. degree in Applied Mathematics from the University of Colorado at Boulder in 1992, and a Ph.D. degree in Applied Mathematics from the University of Colorado, Boulder, supported by a National Science Foundation graduate traineeship under the mentorship of James Meiss, in 1995. His Ph.D. thesis is entitled "Controlling Chaos, Targeting, and Transport" [1].

In 1995, Erik joined the United States Military Academy as an Assistant Professor in the Department of Mathematical Sciences. In 1997, he moved to the United States Naval Academy, where he was an Assistant (1997–2000) and Associate (2000–2002) Professor in the Mathematics Department. In 2002, Erik joined the Mathematics Department at Clarkson University, where he spent the rest of his career. He was promoted to Full Professor in 2006 and was the W. Jon Harrington Professor of Mathematics at the time of his death, with a dual appointment in Electrical and Computer Engineering, and an adjunct appointment in Physics. He was the inaugural director of the Clarkson Center for Complex Systems Science.

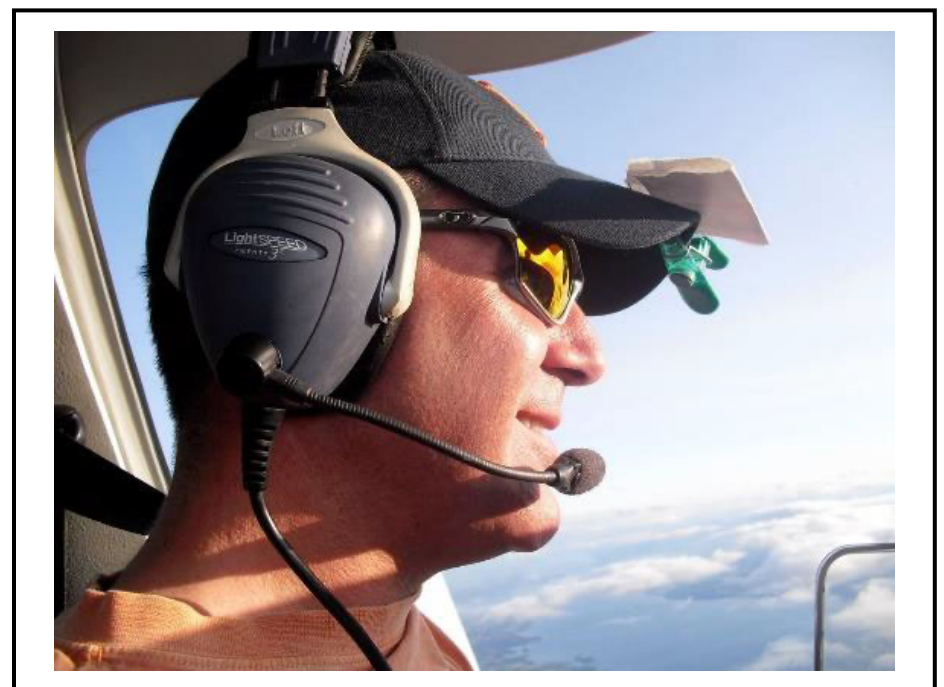

Dr. Erik Bollt, the amateur pilot. From [3].

As Erik described on his website [2], his research interests encompass data-driven analysis of complex systems and dynamical systems, causation inference, machine-learning and data-science methods, stochastic-process modeling and analysis, network science, and numerous applications. Erik published more than 200 research papers, several book chapters, and the book Applied and Computational Measurable Dynamics. During his career, Erik mentored 12 post-doctoral research associates, 16 Ph.D. students, and numerous Master's and undergraduate students. He was astoundingly successful in acquiring extramural funding and was awarded many competitive grants.

In this article, we discuss Erik's areas of research, his service to the community, and describe what Erik was like as a scientist and person.

## Areas of Research

Erik had a wide range of interests in the general area of dynamical and complex systems, bringing rigorous mathematical methods and a keen intuition to the problems he studied, always striving to connect mathematics to real-world applications. He actively sought collaborations with international researchers and had an infectious enthusiasm for new mathematical challenges. His collaborators uniformly comment on his wide range of knowledge and his ability to connect seemingly disparate topics. In this section, we highlight a few of Erik's contributions.

### *Chaos Control and Targeting*

A fundamental question in nonlinear dynamics is whether the rich orbit structure embedded in a chaotic attractor can be exploited for control purposes. In the mid-1990s, Erik developed rigorous methods for targeting, a method of steering a chaotic trajectory towards a desired region of phase space in minimum time by applying only small perturbations [4, 5]. His approach elegantly combined graph-theoretic algorithms with the natural transport properties of chaotic dynamics, showing that the complex web of unstable orbits within a strange attractor can serve as a skeleton for efficient control strategies. What made Erik's work stand out, even in those early days, was his insistence on connecting abstract mathematical ideas to tangible problems: his targeting methods were not simply abstract theoretical results but practical instruments for genuine problem solving.

A particularly striking contribution was his work on the inverse Frobenius–Perron problem: given a desired invariant density, construct a dynamical system, or a perturbation of an existing one, that realizes it [6]. This reversed the classical question of determining the statistical properties of a given map, instead asking how to engineer dynamics with prescribed statistical behavior. Erik showed that this problem admits constructive solutions and connected it to practical questions in chaos control and communication using chaotic signals.

### *Dynamics Mode Decomposition*

In the late 2010's, Erik was instrumental in developing dynamic mode decomposition (DMD), a modern data-driven spectral analysis tool for devising effective models of dynamical systems [7]. At its core, DMD extracts coherent spatiotemporal patterns, referred to as dynamic modes, directly from time-series data, providing a finite-dimensional approximation to the linear Koopman operator even when the true system is nonlinear. Erik was able to seamlessly transition into research in this area because the Koopman operator represents the adjoint to the Frobenius-Perron operator, a topic he had already researched extensively by this point.

The importance of DMD lies in this bridge: it allows researchers to perform spectral analysis, model reduction, prediction, and control using only data, without requiring explicit governing equations. DMD's generalization to extended DMD and later adaptive variants grew from the realization that the choice of observables (the *dictionary*) determines how well one can approximate the infinite-dimensional Koopman operator. Erik's work repeatedly pushed this frontier, clarifying the operator-theoretic foundations, developing practical algorithms, and showing how these spectral decompositions reveal structure in systems ranging from fluid flows to coupled oscillators.

DMD endures as a central technique because it captures the spirit of Erik's approach: take an abstract idea, make it computationally concrete, and use it to illuminate the hidden geometry of complex dynamics.

*Data-Driven Model Discovery*

A core challenge in the analysis of complex systems is to infer causal relationships from observed data without relying on *a priori* models. Erik addressed this problem by introducing causation entropy, a rigorous information-theoretic measure that quantifies the directional causal influence of one time series on another [8–11]. Causation entropy provides a principled criterion for identifying the minimal set of variables that drive a given process, enabling data-driven discovery of causal network structure directly from time-series measurements.

This framework proved to be far more than a theoretical construct. Erik's ideas on causation entropy greatly shaped subsequent work on understanding how coordination and leadership emerge spontaneously in complex human networks. By offering a rigorous, model-free way to detect who is influencing whom in a complex system, his methods opened the door to analyzing motion-capture data from human groups performing synchronization tasks [12], identification of polygenetic factors in gene expression networks [13], and an understanding of how the "asynchronous" functional network of the human brain is affected by alcohol consumption [14]. Erik was genuinely excited by these applications; he saw them as a validation of his conviction that the right mathematical abstraction, rooted in information theory and dynamical systems, could illuminate phenomena far from its original setting, such as understanding something as subtle as human social dynamics.

*Synchronization in Networks*

Erik also brought his characteristic blend of intellectual rigor and appetite for unconventional applications to synchronization, which describes the spontaneous emergence of coordinated behavior of coupled dynamical units. Synchronization is a ubiquitous phenomenon in physics, biology, and engineering. Erik made several contributions to the study of synchronization, where he pushed the research frontiers by analyzing synchronizing units that influenced each other with different weights [15] or through time-dependent and multilayer connections [16]. In Erik's analysis of synchronization on networks, he applied powerful analytical techniques such as master stability functions [17].

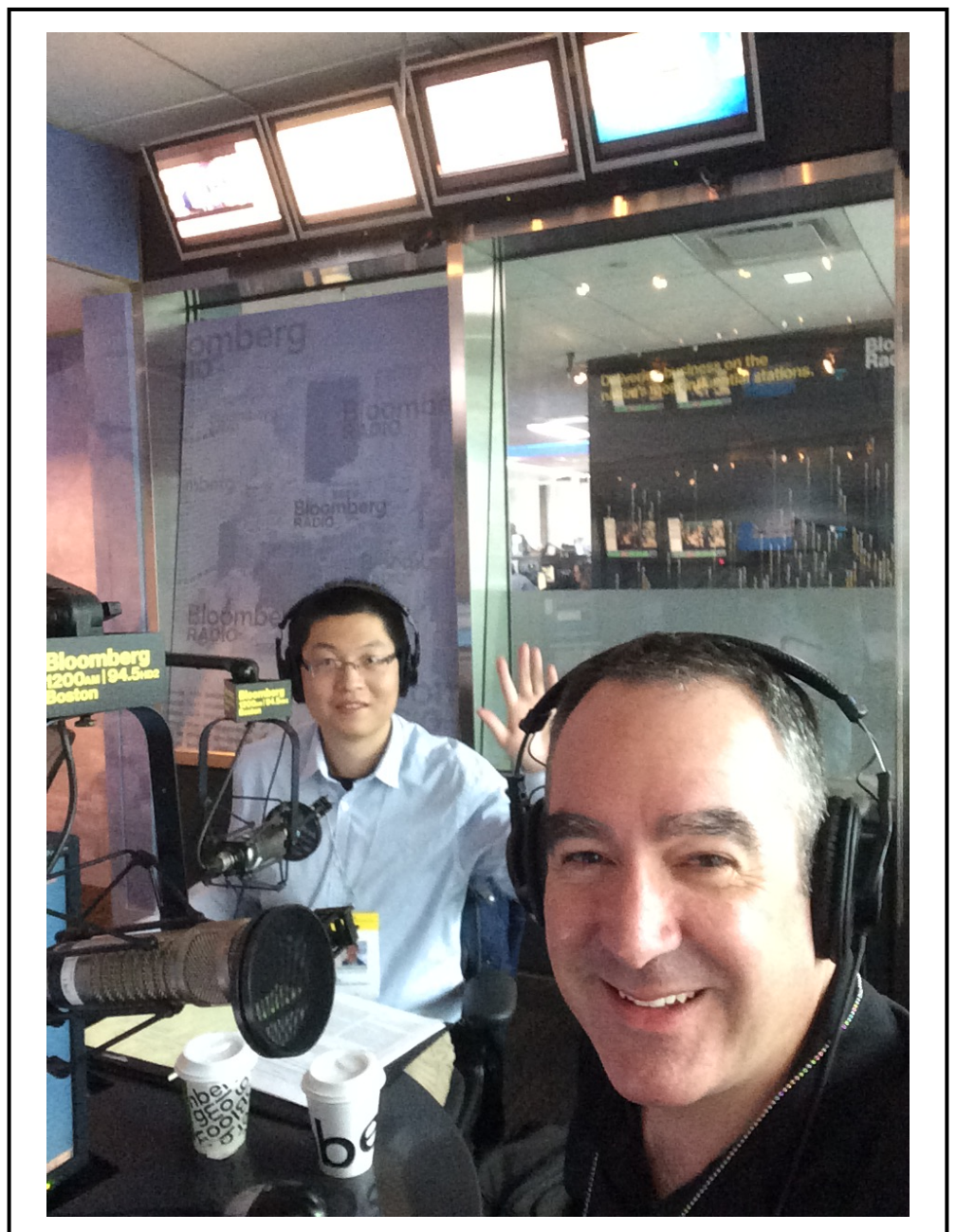

Jie Sun and Erik Bollt during a radio interview on cow synchronization. Photo credit: Erik Bollt.

Synchronization was also the topic of one of Erik's most memorable collaborations: the development of a mathematical model of how cows synchronize their daily activity patterns [18]. The problem was motivated by observations that herds of cows tend to coordinate their lying and standing behavior in striking unison. Erik and his collaborators developed a minimal dynamical-systems model in which each cow's behavioral state (standing or lying) is governed by simple coupling rules. The model takes the form of a piecewise-linear system, which is a type of model that is characterized by abrupt switches in the governing equations. This model provided natural yet highly nontrivial dynamics. The interplay between the switching boundaries and the inter-cow coupling produced rich synchronization patterns that can be analyzed rigorously using tools from the theory of piecewise-smooth dynamical systems. The model demonstrated that synchronization in a cow herd can emerge from purely local interactions, without any need for a leader or global signal.

What made the project quintessentially Erik was the joyful way that he approached the collaboration. During the weeks of daily working sessions, the team would sit together in front of a computer, build and test the model in real time, and try to imagine what it was like to be a cow. Erik relished getting into the fine details of the problem, questioning every modeling assumption and pushing for the simplest formulation that still captured the essential dynamics. The result was a paper [18] that resonated far beyond the dynamical-systems community. A few years later, it was even featured on the Bloomberg 'Odd Lots' podcast, where Erik and the team discussed with infectious enthusiasm what herding behavior in cows might teach us about coordination in financial markets [19].

*Machine Learning*

Erik's intellectual curiosity naturally led him to study the burgeoning intersection of machine learning and complex nonlinear dynamical systems. One work developed and studied Kolmogorov-Arnold networks (KANs), which addresses [20] the limitations of sparse optimization for data-driven model discovery, specifically, the assumption that governing equations must contain only a small number of elementary functions. The goal was to develop a systematic framework for identifying elementary function bases, which allows for the interpretation of the network processing capabilities. This is opposed to deep neural networks, which can provide accurate predictions, but are opaque black boxes.

In the collaborative project, Erik helped to define a KAN-based framework that broke free from the sparsity constraint. This approach allowed the network to accurately capture complex behaviors while retaining the interpretability that traditional neural networks lack. Erik was particularly excited by the insight this provided into the underlying dynamics, insights usually lost in black-box approximations. This work demonstrated an interesting non-uniqueness, revealing that many approximate models can generate the same invariant sets that matched the underlying statistics of the original dynamical system, such as Lyapunov exponents. Erik helped make a connection of this result to previous work on shadowing numerical trajectories in chaotic systems, linking this modern machine-learning approach back to the fundamental chaos theory he loved.

Erik also studied a different machine learning model known as reservoir computing, which is well suited for predicting the dynamics of a system using only observed time series data. Here, the artificial neural network consists of an input layer, a hidden recurrent unit, and an output layer. The neurons in the hidden layer are nonlinear dynamical systems that naturally have memory. To make the problem tractable, the weights of the input and hidden unit are not trained: they are assigned randomly and kept fixed. Only the output layer is trained using the observed data, which is input to the network. The output layer is a linear superposition of the internal states, which requires linear optimization rather than the complex nonlinear optimization needed for other machine learning algorithms. This convex optimization problem is solved using regularized linear regression for which there are efficient numerical methods.

Erik liked the simplicity of the model and the fact that it requires only a short training dataset. But he was bothered by the random assignment of the weights on the input layer and the hidden unit. To address this problem, Erik made the neurons linear and found that the linear artificial neural network can perform very well for a short forecasting horizon [21]. Furthermore,

Erik showed that the functions of the network can be separated: the memory can be realized using a delay line on the input, and the nonlinearity can be moved to the output layer.

Taking this line of research a step further, Erik and collaborators showed that the input delay line can be shallow, only including data from a small, finite set of steps in the past, and the output layer is a linear superposition of functionals [22]. He showed accurate forecasting of chaotic dynamical systems can be obtained with a neural network much smaller than previous work and required fewer observations for model training.

*Symbolic Dynamics*

One of Erik's most enduring legacies lies in his transformative work on symbolic dynamics, particularly the challenge of assigning symbols to chaotic attractors. This task had long been hindered by the intricate, fractal nature of generating partitions in non-uniformly hyperbolic systems.

To tackle this problem, Erik realized that the unstable periodic orbits (UPOs) that form the skeleton of a chaotic attractor can be used to progressively approximate the partition itself [23]. His fundamental observation was that the coarse topological features of chaotic attractors are typically revealed by a relatively small number of short UPOs, while increasingly longer orbits serve to refine features without altering the general structure. Based on this, he pioneered an efficient algorithm that uses proximity functions in phase space to assign symbols to orbit points. By exploiting the principle that points on longer UPOs are likely to share symbolic assignments with nearby points on shorter orbits, he enabled the computation of generating partitions for systems far more complex than one-dimensional maps, such as the two-dimensional Ikeda map. This work stands as a benchmark for generating partitions in the study of chaotic symbolic dynamics.

In another foundational contribution, Erik brought his characteristic rigor to the evaluation of threshold-crossing analysis. At a time when experimentalists frequently encoded chaotic time-series using arbitrary partitions, Erik recognized the mathematical perils of this approach. He argued that the generating partition of a chaotic system connects primary tangencies and cannot be approximated by simple lines used in threshold crossings. Through rigorous analysis of the tent map, he demonstrated that *sample-path* symbolic dynamics derived from arbitrary partitions could severely misrepresent the dynamical system, leading to diminished topological entropy and non-uniqueness [24]. In a surprising and deep mathematical finding, he showed that the topological entropy, viewed as a function of partition misplacement, behaved in a non-monotone, devil's staircase-like manner. These insights provided a critical quantitative framework for the field, insisting on mathematical validity in the face of popular heuristic techniques.

**Organizing the Community**

Erik Bollt was an active and inspiring member of the applied nonlinear dynamics and complex systems communities. He was a key presence at the biennial SIAM Conference on Applications of Dynamical Systems (*i.e.*, "the Snowbird meeting"), which was held in the spring of odd years. During the SIAM meeting, Erik was a familiar presence with conversations often turning to possible applications of his fundamental research topics. Erik had an infectious way of sketching connections on napkins and whiteboards, always eager to see how an abstract result might find a home in a new context. Those discussions captured something essential about him: a deep mathematical seriousness paired with a playful curiosity.

Erik was also a major contributor to an annual series of complex-systems conferences in the northeastern United States, which played a key role in developing the local complex-systems community. This contribution was formalized by his founding in 2017 of $C^2S^3$: The Clarkson Center for Complex Systems Science. The center was instrumental in developing interdisciplinary research in dynamical systems across the campus and hosted many international visitors and sponsored local workshops.

Mario di Bernardo, Mattia Frasca, and Erik Bollt at the 2023 SIAM Conference on Applied Dynamical Systems, Portland, OR, USA. Photo credit: Mario di Bernardo.

**Erik Bollt the Person**

Erik possessed a rare combination of deep mathematical insight, genuine enthusiasm for collaboration, and generosity in sharing ideas. His ability to bridge theory and application, his openness to interdisciplinary approaches, and his gift for seeing unexpected connections made him an inspiration to colleagues worldwide. In addition to being a brilliant mathematician and scientist, Erik was also a great collaborator, and a cherished mentor to his students and postdocs. Erik exuded a rational positivity that brought light to his mentees, and he brought this same positivity home with him. As a father, he supported his children's varied interests by tagging along with each of them on their individual journeys and helping whenever and however he could. His pride in his children and their accomplishments shone through every interaction. Erik: a father, mentor, scientist, athlete, pilot, and even a video-game developer [25].

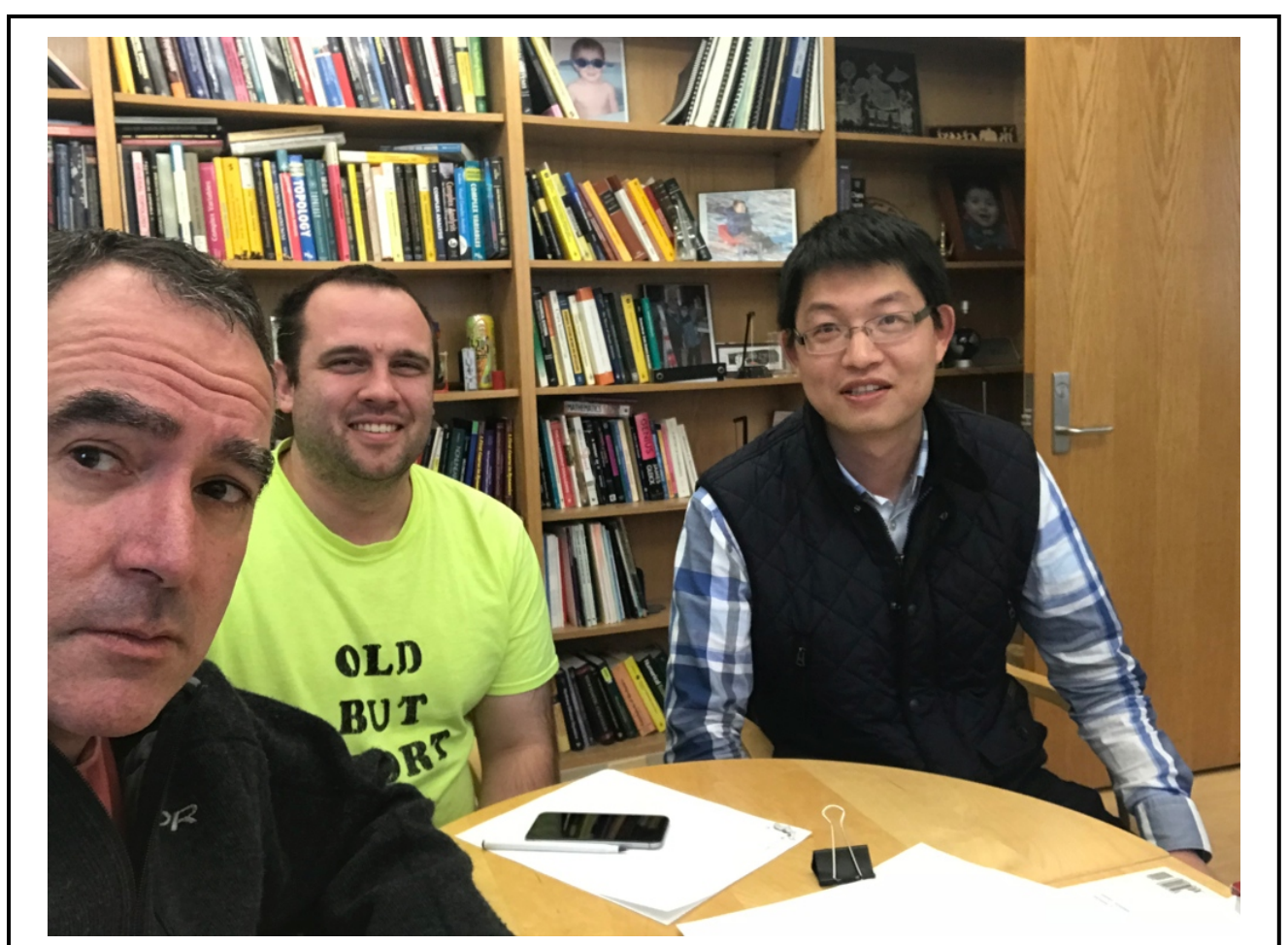

Erik Bollt, Jeremie Fish, and Jie Sun in Erik's Clarkson office. Photo credit: Erik Bollt